\documentclass[12pt]{amsart}
\usepackage{amscd,amssymb}
\usepackage[graph,frame,poly,arc]{xy}  
\usepackage[plainpages,backref,urlcolor=blue]{hyperref}

\theoremstyle{plain}
\newtheorem{thm}[subsection]{Theorem}
\newtheorem{lem}[subsection]{Lemma}
\newtheorem{prop}[subsection]{Proposition}
\newtheorem{cor}[subsection]{Corollary}

\theoremstyle{definition}
\newtheorem{rk}[subsection]{Remark}
\newtheorem{definition}[subsection]{Definition}

\numberwithin{equation}{section}
\newcommand{\F}{{\mathcal F}}
\newcommand{\A}{{\mathcal A}}
\newcommand{\B}{{\mathcal B}}

\newcommand{\Z}{\mathbb{Z}}

\newcommand{\C}{\mathbb{C}}
\newcommand{\PP}{\mathbb{P}}

\DeclareMathOperator{\rank}{rank}

\begin{document}

\title [On the exponents of free and nearly free projective plane curves]
{ On the exponents of free and nearly free projective plane curves}

\author[Alexandru Dimca]{Alexandru Dimca$^1$}
\address{Univ. Nice Sophia Antipolis, CNRS,  LJAD, UMR 7351, 06100 Nice, France. }
\email{dimca@unice.fr}

\author[Gabriel Sticlaru]{Gabriel Sticlaru}
\address{Faculty of Mathematics and Informatics,
Ovidius University,
Bd. Mamaia 124, 900527 Constanta,
Romania}
\email{gabrielsticlaru@yahoo.com }

\thanks{$^1$ Partially supported by Institut Universitaire de France.}



\begin{abstract} We show that all the possible pairs of integers occur as exponents for free or nearly free irreducible plane curves and line arrangements, by producing only two types of simple families of examples. The topology of the complements of these curves and line arrangements is also discussed, and many of them are shown not to be $K(\pi,1)$ spaces.

\end{abstract}
 
\maketitle


\section{Introduction} 

Let $C:f=0$ be a reduced free curve of degree $d$ with exponents $d_1 \leq d_2$. Then one has
$$d_1+d_2=d-1,$$
see for instace \cite{DS14}, \cite{DStFD}.
We assume that $C$ is not a union of lines passing through one point, which is equivalent to asking $d_1>0$. When $C$ is irreducible, one also knows that $d_1 \geq 2$, see \cite{DStFD}.
The following natural question seems not to have a clear answer in the existing literature to the best of our knowledge: given a degree $d$, which are the pairs $(d_1,d_2)$ of exponents which may occur for a free curve of degree $d$ ?
If we fix the minimal degree $d_1$ of a syzygy involving the partial derivatives $f_x,f_y,f_z$ of $f$,
then the freeness of the curve $C$ can be characterized by a maximality property of the corresponding total Tjurina number $\tau(C)$, see Theorem 3.2 in du Plessis and Wall's paper \cite{duPCTC}, and its discussion in \cite{Dmax}. However, in these papers it is not clear in our opinion for which values of $d_1$ this maximality property can hold.

To our surprise, especially in view of the recent paper \cite{B+} showing that the free and nearly free curves can be quite wild, in particular can fail to be rational, all possible values of the exponents can be realized just by looking at very restrictive classes of free curves. Here are the results, first for irreducible curves, and then for line arrangements.

\begin{thm}
\label{thm2} Consider the irreducible curve $C : f=0$ of degree $d$  in $\PP^2$ given by
$$f(x,y,z)= x^d+x^{d_1}y^{d_2+1}+y^{d-1}z,$$
for $2 \leq d_1 <d/2$ and $d_2=d-1-d_1$. Then the following holds.

\begin{enumerate}

\item The curve $C$ is free with exponents $(d_1,d_2)$.

\item The curve $C$ is rational and has a unique cusp, which has a unique Puiseux pair, namely $(d-1,d)$.

\end{enumerate}

\end{thm}

 For two integers $i \leq j$ we define a homogeneous polynomial in $\C[u,v]$ of degree $j-i+1$ by the formula 
\begin{equation}
\label{gij}
 g_{i,j}(u,v)=(u-iv)(u-(i+1)v) \cdots (u-jv).
\end{equation}

\begin{thm}
\label{thm1} Consider the line arrangement $\A : f=0$ of $d\geq 3$ lines in $\PP^2$ given by
$$f(x,y,z)= xg_{1,d_1}(x,y)g_{1,d_2}(x,z)=0,$$
for $1 \leq d_1 <d/2$ and $d_2=d-1-d_1$. Then the following holds.

\begin{enumerate}

\item The line arrangement $\A$ is free with exponents $(d_1,d_2)$.

\item The line arrangement $\A$ has at most two points of multiplicity $>2$.
\end{enumerate}

\end{thm}

The corresponding results for nearly free irreducible curves and line arrangements are stated in 
Theorems \ref{thm3} and \ref{thm4}. In Proposition \ref{prop5} we show, using the results by 
du Plessis and Wall's paper \cite{duPCTC2} on plane curves with symmetry, that for nearly free line arrangements one has $d_1\geq 2$, i.e. the same condition as for the irreducible free curves.

The proofs in the irreducible case rely on the verification of Saito's criterion of freeness and of its extension to nearly free curves established in \cite{Dmax}. In the case of line arrangements, we can use either the geometric results by Faenzi and Vall\` es in \cite{FV} and by Abe \cite{A}, or directly the characterization of (nearly) freeness in terms of the maximality of the total Tjurina number from
\cite{duPCTC} (in the free case) and \cite{Dmax} (in the nearly free case).

\medskip

In the last section we discuss the topology of the complements in $\PP^2$ of the families of curves introduced in Theorems \ref{thm2}, \ref{thm1}, \ref{thm3} and \ref{thm4}. We show that most of these complements are not $K(\pi,1)$ spaces.
We thank Mike Falk for his help in proving Proposition \ref{propANF}.

\section{Proofs of Theorems \ref{thm2} and \ref{thm1}}

\subsection{Proof of Theorem \ref{thm2}}

One can easily find the following two syzygies involving the partial derivatives of $f$.

\begin{equation}
\label{s1}
y^{d_1}f_y-((d_2+1)x^{d_1}+(d-1)y^{d_1-1}z)f_z=0,
\end{equation}

and

\begin{equation}
\label{s2}
(d_2+1)^2y^{d_2}f_x-d ((d_2+1)x^{d_2}-(d-1)x^{d_2-d_1}y^{d_1-1})f_y-  
\end{equation}
$$-(d_1(d_2+1)^2x^{d_1-1}y^{d_2-d_1+1}+d(d-1)^2x^{d_2-d_1}y^{d_1-2}z^2)f_z=0.$$

We can conclude that $C$ is a free curve in two ways. The simpliest way is to use Lemma 1.1 in \cite{ST}, which says that a curve of degree $d$ having two independent syzygies of degree $d_1$ and $d_2$ such that $d_1+d_2 \leq d-1$ is free with exponents $(d_1,d_2)$.

Alternatively, one may use Saito's criterion, see for instance  \cite{KS} or \cite{Yo}, namely compute the determinant $\Delta$ of the $3 \times 3$ matrix obtain using for the first line $x,y,z$, for the second line the coefficients of the syzygy \eqref{s1} and for the third line the coefficients of the syzygy \eqref{s2}. One obtains that
$\Delta=c \cdot f$, for a nonzero constant $c$, which again implies that $C$ is free with exponents $(d_1,d_2)$.

\subsection{Proof of Theorem \ref{thm1}}  Note that the line arrangement $\A$ has the following two points of multiplicity (possibly) strictly greater than 2: the point $A=(0:0:1)$ of multiplicity
$d_1+1$ and the point $B=(0:1:0)$ of multiplicity $d_2+1$. 
Then we apply Theorem 2 in \cite{FV} by taking $k=d_1$ and $r=d_2-d_1$. We clearly have $d=2k+r+1$ and the line arrangement $\A$ has a point of multiplicity $h=d_1+1 \in [k,k+r+1]$.
It follows that $\A$ is free with exponents $(d_1,d_2)$ if and only if the second Chern class
$c_2(T_{\PP^2}(-log \A))$ is given by $d_1d_2$. Note that the vector bundle $T_{\PP^2}(-log \A)$ corresponds to the vector bundle denoted by $T<\A>(-1)$ in \cite{DS14}, and hence the previous claim is compatible with the formula
$$c_2(T<\A>(-1))=(d-1)^2-\tau(\A),$$
given in \cite{DS14}, see for instance the formula (3.2).
The formulas (4) and (5) in Remark 2.2 in \cite{FV} imply that we have the following formula

\begin{equation}
\label{c2}
c_2(T_{\PP^2}(-log \A))= \sum_{j\geq 2}(j-1)n_j -d+1,
\end{equation}
where $n_j$ denotes the number of points of multiplicity $j$ in the line arrangement $\A$.
As we have seen above, there is one point of multiplicity $d_1+1$, one point of multiplicity $d_2+1$ and $d_1d_2$ additional nodes. It follows that
$$c_2(T_{\PP^2}(-log \A))=d_1d_2+d_1+d_2-d+1=d_1d_2$$
and hence $\A$ is free with exponents $(d_1,d_2)$.

\section{Exponents of nearly free plane curves}  

In this section we discuss the exponents of nearly free curves, a notion introduced in \cite{DStNF}
in our study of rational cuspidal curves. Corollary \ref{corNF} given below shows once more that it is natural to treat the free curves and the nearly free curves together.

\begin{thm}
\label{thm3} 

Consider the irreducible curve $C : f=0$ of degree $d$  in $\PP^2$ given  by
$$f(x,y,z)= x^d+x^{d_2+1}y^{d_1-1}+y^{d-1}z,$$
for $1 \leq d_1 \leq d/2$ and $d_2=d-d_1$. Then the following holds.

\begin{enumerate}

\item The curve $C$ is nearly free with exponents $(d_1,d_2)$.

\item The curve $C$ is rational and has a unique cusp, which has a unique Puiseux pair, namely $(d-1,d)$.

\end{enumerate}

\end{thm}

\proof The case $d_1=1$ follows from Proposition 5.1 in  \cite{DStNF}. So from now on we suppose $d_1 \geq 2$.
Then it is easy to find the following syzygies, where we set $k=d_1-1\geq 1$.

\begin{equation}
\label{s11}
kxy^{k}f_x-(dx^ky+(d-k)y^{k+1})f_y+(d(d-1)x^kz+(d-1)(d-k)y^kz)f_z=0.
\end{equation}

\begin{equation}
\label{s12}
-ky^{d_2}f_x+dx^{k-1}y^{d-2k}f_y+(k(d-k)x^{d_2}-d(d-1)x^{k-1}y^{d-2k-1}z)f_z=0.
\end{equation}

The last syzygy $af_x+bf_y+cf_z=0$ is more difficult to determine. Such a syzygy implies that $a$ is divisible by $y^{k-1}$, say $a=y^{k-1}a_1$. It follows that
$$x^{d-k-1}[a_1(dx^k+(d-k)y^k)+kbx]=-y^{d-k-1}((d-1)bz+cy).$$
Hence there is a linear form $L$ such that $a_1(dx^k+(d-k)y^k)+kbx=y^{d-k-1}L$ and
$(d-1)bz+cy=-x^{d-k-1}L$. By taking $a_1=qx^{d-2k}+y^{d-2k-1}z$, with $q$ to be determined, one gets the following solution.

\begin{equation}
\label{s13}
 a=y^{k-1}(q x^{d-2k}+y^{d-2k-1}z),
\end{equation}
\begin{equation}
\label{s14}
 b=-[qdx^{d-k-1}+q(d-k)x^{d-2k-1}y^k+dx^{k-1}y^{d-2k-1}z]/k,
\end{equation}
\begin{equation}
\label{s15}
 c=(d-1)[q(d-k)x^{d-2k-1}y^{k-1}z+dx^{k-1}y^{d-2k-2}z^2]/k,
\end{equation}
where
$$q=\frac{(d-k)k}{d(d-1)}.$$

Let $r_1$ (resp. $r_2$) be the vector constructed using the coefficients of $f_x,f_y,f_z$ in the syzygy \eqref{s11} (resp. \eqref{s12}).
Let  $r_3$ be the vector $(a,b,c)$, with $a,b,c$ defined above.
One can then verify the relation
$$-qx^{d_2-d_1+1}r_1+zr_2+kyr_3=0,$$
which implies that $C$ is nearly free with exponents $(d_1,d_2)$ using the characterization of nearly free curves in Theorem 4.1 in \cite{Dmax}.

\endproof

\begin{rk}
\label{special} 
Using the same approach as in the proof above one can show that the curve $C: f=x^{2k}+x^ky^k+y^{2k-1}z=0$, resp. $C':f'=x^{2k+1}+x^{k+1}y^k+y^{2k}z=0$,
is also nearly free with exponents $(k,k)$, resp. $(k,k+1)$. Note that the syzygies given in the proof of Theorem \ref{thm3} do not apply to these two special cases, since some exponents there become negative.

\end{rk}

This remark and Theorems  \ref{thm2}, \ref{thm3} imply the following.

\begin{cor}
\label{corNF}
Consider the rational cuspidal curve $C:f=x^d+x^ay^{b}+y^{d-1}z=0$ with $a\geq 2$, $b\geq 1$.
 and $a+b=d$. Then the curve $C$ is free for $a<b$ with exponents $(a,b-1)$ and nearly free for $a \geq b$. The exponents in the latter case are $(b+1,a-1)$ if $a \geq b+2$ and $(b,a)$ if $b \leq a \leq b+1$.
\end{cor}

In particular, for $d=2d'$ even, the pairs $(a,b)=(d',d')$ and $(a,b)=(d'+1,d'-1)$ lead to the same exponent $(d',d')$, while for $d=2d'+1$ odd, the pairs $(a,b)=(d'+1,d')$ and $(a,b)=(d'+2,d'-1)$
lead to the same exponent $(d',d'+1)$.

\bigskip

Now we consider the nearly free line arrangements. Essentially we use the same arrangement as in Theorem \ref{thm1}, but we add a new line and create in this way a controlled number of triple points.

\begin{thm}
\label{thm4} Consider the line arrangement $\A : f=0$ of $d=d_1+d_2 \geq 4$ lines in $\PP^2$ given by
$$f(x,y,z)= x(y-z)g_{1,d_1-1}(x,y)g_{2,d_2}(x,z)=0$$
for $2 \leq d_1 \leq d_2$. Then the following holds.

\begin{enumerate}

\item The line arrangement $\A$ is nearly free with exponents $(d_1,d_2)$.

\item The line arrangement $\A$ has at most two points of multiplicity $>3$.
\end{enumerate}

\end{thm}

\proof As in the proof of Theorem \ref{thm1}, and using the same notation, it follows that we have
\begin{equation}
\label{tau}
 (d-1)^2-\tau(\A)=\sum_{j\geq 2}(j-1)n_j -d+1.
\end{equation}
Let us consider the multiple points of the arrangement $\A$. The point $A=(0:0:1)$ has multiplicity $d_1$, the point $B=(0:1:0)$ has multiplicity $d_2$. There are $d_1-2$ triple points on the line $y-z=0$, namely the points $(2:1:1),...,(d_1-1:1:1)$ and $d-1-2(d_1-2)=d_2-d_1+3$ double points.

And the lines passing through $A$ and $B$ intersect in $(d_1-1)(d_2-1)-(d_1-2)$ additional nodes.
Hence
$$
 \sum_{j\geq 2}(j-1)n_j -d+1=$$
$$= (d_1-1)+(d_2-1)+[(d_1-1)(d_2-1)-(d_1-2)+d_2-d_1+3]+2(d_1-2)-d+1=$$
$$d_1(d_2-1)+1,$$
which implies 
\begin{equation}
\label{tau1}
 \tau(\A)=(d-1)^2-d_1(d_2-1)-1.
\end{equation}
By Theorem 3.1 in \cite{Dmax}, in order  to prove that $\A$ is nearly free it is enough to show that $d_1$ is the minimal degree of a nontrivial syzygy
$$R: ~ ~     ~ ~  ~  ~  af_x+bf_y+cf_z=0.$$
In other words we have to show that

\noindent (i) there is  a syzygy $R$ with $a,b, c$ homogeneous of degree $d_1$, and

\noindent (ii) any such syzygy $R$ with $a,b, c$ homogeneous of degree $m$ satisfies $m \geq d_1$.

We start with one. Note that one can write
$$f_y=-f \Big ( \frac{1}{x-y}+\frac{2}{x-2y}+ \cdots +\frac{d_1-1}{x-(d_1-1)y}-\frac{1}{y-z}\Big ).$$
This implies that 
$$f_y=f\frac{P}{Q}$$
where $P \in \C[x,y,z]$ is homogeneous of degree $d_1-1$,  $Q \in \C[x,y,z]$ is homogeneous of degree $d_1$ and $P$ and $Q$ have no common factor. It follows that
$dQf_y=dfP=xPf_x+yPf_y+zPf_z$, or equivalently
$$R_1: ~ ~     ~ ~  ~  ~xPf_x+(yP-dQ)f_y+zPf_z=0.$$
Note that $a_1=xP,$ $ b_1=yP-dQ$ and $c_1=zP$ have no common factor in $\C[x,y,z]$, and $a_1,b_1,c_1$ are homogeneous of degree $d_1$, hence this syzygy $R_1$ satisfies (i).

Suppose now that there is a non trivial syzygy $R$ of degree $m$ with $m \leq d_1-1$. Since $R_1$ cannot be a multiple of $R$, it follows that we can apply Lemma 1.1 in \cite{ST} or the vanishing (4.1) in \cite{Dmax}, and deduce that $m=d-d_1-1=d_2-1 \leq d_1$ and $\A$ is free with exponents $(d_2-1,d_1)$. But this is impossible, since a free curve $C$ with exponents $(d_2-1,d_1)$
has a Tjurina number $\tau(C)= (d-1)^2-d_1(d_2-1)$, see for instance \cite{DStFD}. But this equality  fails for $\A$, hence the existence of a syzygy $R$ of degree $m \leq d_1-1$ is impossible.
This completes the proof of Theorem  \ref{thm4}.

\endproof

The following result completes the picture.

\begin{prop}
\label{prop5} There is no nearly free line arrangement with exponents $(d_1,d_2)$ if $d_1=1.$

\end{prop}

\proof

Suppose $\A$ is such an arrangement. Then as in formula \eqref{tau1}, one should have
$$\tau(\A)=(d-1)^2-(d-2)-1=d^2-3d+2.$$
Proposition 3.1 in \cite{duPCTC2} implies that $\A$ is invariant by a semisimple subgroup $H$ in $Gl_3(\C)$. 
In fact, one can assume that $d =d_1+d_2 \geq 4$, since all line arrangements with $d<4$ are clearly free. Then Proposition 1.2 in \cite{duPCTC2} implies that $\dim H =1$. More precisely, one can take
$H=\C^*$ and the action of $H$ on $\PP^2$ is given by
$$t\cdot (x:y:z)=(x:t^by:t^cz),$$
for some integers $b$ and $c$, not both of them zero, see the beginning of section 2 in \cite{duPCTC2}.

Note that a line arrangement is invariant by such an action if and only if any line in the arrangement is invariant. It is clear that the lines $x=0$, $y=0$ and $z=0$ are invariant under this diagonal action. But the line arrangement $xyz=0$ is free.

A line $L: px+qy=0$ with both $p$ and $q$ nonzero is invariant by $H$ if and only if 
 $b=0$. It follows that such a group action leaves invariant a line arrangement consisting of some  lines through the point $A=(0:0:1)$ and possibly the line $z=0$. But these two types of arrangements are easily seen to be free.

Finally, a line $L:px+qy+rz=0$ with all $p,q,r$ nonzero is never invariant by a nontrivial action of $H=\C^*$ as above.

This completes the proof of Proposition \ref{prop5}.

\endproof

\section{The topology of the complements}  

In this section we look at the topology of the complements $U=\PP^2 \setminus C$, where $C$ is one of the free or nearly free curves discussed above. The case of the irreducible curves is covered by the following slightly more general result.

\begin{prop}
\label{propIRR} Let $C:f=0$   be a rational cuspidal curve of degree $d>1$ in $\PP^2$, having a point of multiplicity $d-1$ and set  $U=\PP^2 \setminus C$. Then the following hold.
\begin{enumerate}

\item $\pi_1(U)=H_1(U,\Z)=\Z/d\Z.$

\item $H_2(U,\Z)=0.$

\item The Milnor fiber $F$ associated to the plane curve $C$ is homotopically equivalent to a bouquet of $(d-1)$ spheres $S^2$.
\end{enumerate}
In particular, $U$ is not a $K(\pi,1)$-space.

\end{prop}

\proof
The first claim follows from the fact that the existence of a point of multiplicity $(d-1)$ forces the fundamental group $\pi_1(U)$ to be abelian, see for instance  \cite[Corollary (4.3.8), p.124]{D1}.

The second claim follows from the fact that $H_2(U,\Z)$ is torsion free (since the variety $U$, being affine, has the homoopy type of a CW-complex of dimension 2) and moreover, one has for the Euler number of $U$:
$$E(U)=E(\PP^2)-E(C)=3-2=1,$$
since $C$ is rational cuspidal, hence homeomorphic to $\PP^1$.

To prove the last claim, note that the Milnor fiber $F$ is given by $f(x,y,z)-1=0$ in $\C^3$ and it is the cyclic $d$-fold cover of $U$. It follows from (1) that $F$ is simply connected. Since $F$ is an affine surface, it follows that it has has the homoopy type of a CW-complex of dimension 2, hence $F$ must be a bouquet of spheres $S^2$. Moreover, one has $E(F)=dE(U)=d$, which implies that there are $(d-1)$ spheres in this bouquet.

Finally, one has $\pi_2(U)=\pi_2(F)=H_2(F,\Z)=\Z^{d-1} \ne 0$, and hence $U$ is not a $K(\pi,1)$-space.
\endproof

Now we consider the complements of the line arrangements described in Theorems \ref{thm1} and \ref{thm4}.

\begin{prop}
\label{propAF} Consider the line arrangement $\A : f=0$  in $\PP^2$ given by
$$f(x,y,z)= xg_{1,d_1}(x,y)g_{1,d_2}(x,z)=0,$$
for $1 \leq d_1 \leq d_2$. If $M(\A)$ denotes the complement of $\A$ in $\PP^2$, then one has
$$M(\A)=(\C \setminus \{ d_1 \text { points} \}) \times (\C \setminus \{ d_2 \text { points} \}).$$
In particular,  $M(\A)$ is a $K(\pi,1)$-space where $\pi= \F_{d_1} \times \F_{d_2}$, with
$\F_m$ the free group on $m$ generators.

\end{prop}

\proof We delete first the line $x=0$, and we get the affine plane with coordinates $y,z$.
The trace $\B$ of $\A$ on this affine plane is given by the equation $g_{1,d_1}(1,y)g_{1,d_2}(1,z)=0,$
which implies our claim.

\endproof

\begin{prop}
\label{propANF} Consider the line arrangement $\A : f=0$  in $\PP^2$ given by
$$f(x,y,z)= x(y-z)g_{1,d_1-1}(x,y)g_{2,d_2}(x,z)=0$$
for $2 \leq d_1 \leq d_2$. If $M(\A)$ denotes the complement of $\A$ in $\PP^2$, then  $M(\A)$ is not a $K(\pi,1)$-space if either $d_1=2$ or $d_1=d_2=3$.

\end{prop}

\proof

Consider first the case $d_1=2$. We proceed as above, i.e. we delete first the line $x=0$, and we get the affine plane with coordinates $y,z$.
The trace $\B$ of $\A$ on this affine plane is given by the equation 
$$(y-z)(y-1)(z-\frac{1}{2})(z-\frac{1}{3})\cdots (z-\frac{1}{d_2})=0.$$
 This affine arrangement $\B$ is the complexification of a real arrangement having a triangular chamber with double point vertices, and hence by \cite[Example 2.7]{F}  $M(\B)$, the complement of $\B$ in $\C^2$, is not  a $K(\pi,1)$-space. Since $M(\A)=M(\A)$, this proves our claim in the first case.

Assume now that $d_1=d_2=3$.
Then the trace $\B$ of $\A$ on the affine plane is given by the equation 
$$(y-z)(y-1)(y-\frac{1}{2})(z-\frac{1}{2})(z-\frac{1}{3})=0.$$
This affine arrangement is  lattice isotopic to the arrangement
$$\B': (y+z-1)(y-1)(y-\frac{1}{2})(z-\frac{1}{2})(z-\frac{1}{3})=0,$$
see \cite[Definition 5.27]{OT} and use the family of lines 
$$(y-\frac{1}{2})-e^{\pi i t}(z-\frac{1}{2})=0,$$
for $t \in [0,1].$
 In particular their complements are homeomorphic, see \cite{Ra}.
 The arrangement $\B'$ is the complexification of a real arrangement having a triangular chamber with double point vertices, and hence as above by \cite[Example 2.7]{F}  $M(\B')$ is not  a $K(\pi,1)$-space. A similar arrangement is discussed in \cite{FR}, just after Theorem 2.4 and in \cite[Example 6.9]{ACM}, where it is mentioned that in fact $\pi_1(M(\B'))$ is the Stallings group.

\endproof


\begin{thebibliography}{00}



\bibitem{A} T. Abe, Roots of characteristic polynomials and intersection 
points of line arrangements. J. of Singularities 8 (2014),100-117.

\bibitem{ACM}   E. Artal Bartolo, J.I. Cogolludo-Agust\'in, D. Matei, Quasi-projectivity, Artin-Tits groups, and pencil maps. Topology of algebraic varieties and singularities, 113--136, Contemp. Math., 538, Amer. Math. Soc., Providence, RI, 2011.

\bibitem{B+} E. Artal Bartolo, L. Gorrochategui, I. Luengo, A. Melle-Hern\' andez,
On some conjectures about free and nearly free divisors, arXiv:1511.09254






\bibitem{D1} A.~Dimca, {\em Singularities and topology of hypersurfaces}, Universitext, Springer-Verlag, New York, 1992. 



\bibitem{Dmax} A. Dimca, Freeness versus maximal global Tjurina number for plane curves, arXiv:1508.0495.










\bibitem{DS14} A. Dimca, E. Sernesi,  Syzygies and logarithmic vector fields along plane curves,
Journal de l'\'Ecole polytechnique-Math\'ematiques 1(2014), 247-267.





\bibitem{DStFD} A. Dimca, G. Sticlaru, Free divisors and rational cuspidal plane curves, arXiv:1504.01242v4.

\bibitem{DStNF} A. Dimca, G. Sticlaru, Nearly free divisors and rational cuspidal curves, arXiv:1505.00666v3.




\bibitem{duPCTC} A.A. du Pleseis,  C.T.C. Wall, Application of the theory of the
discriminant to highly singular plane curves, Math. Proc. Camb.
Phil. Soc.,  126(1999), 259-266. 

\bibitem{duPCTC2} A. A. du Plessis and C. T. C.Wall, Curves in $P^2(\C)$ with 1-dimensional symmetry,
Revista Mat Complutense 12 (1999), 117--132.


\bibitem{FV}  D. Faenzi, J. Vall\`es, 
{Logarithmic bundles and line arrangements, an approach via the standard construction}, J. London.Math.Soc.
{90} (2014) {675--694}.

\bibitem{F} M. Falk, $ K(\pi,1)$ arrangements. Topology 34 (1995), no. 1, 141--154.

\bibitem{FR}  M. Falk, R. Randell,  On the homotopy theory of arrangements. II. Arrangements--Tokyo 1998, 93--125, Adv. Stud. Pure Math., 27, Kinokuniya, Tokyo, 2000.

\bibitem{OT} P. Orlik and H. Terao, {\em Arrangements of Hyperplanes,}Springer-Verlag, Berlin Heidelberg New York, 1992.

\bibitem{Ra}  R. Randell, Lattice-isotopic arrangements are topologically isomorphic, Proc. Amer. Math. Soc. 107(1989), 555-559.

\bibitem{KS} K. Saito, Theory of logarithmic differential forms and logarithmic vector fields, J. Fac. Sci. Univ. Tokyo Sect. IA Math. 27 (1980), no. 2, 265-291.









\bibitem{ST} A. Simis, S.O. Tohaneanu, Homology of homogeneous divisors, Israel J. Math. 200 (2014), 449-487.


\bibitem{Yo} M. Yoshinaga, Freeness of hyperplane arrangements and related topics, Annales de la Facult\'e des Sciences de Toulouse, vol. 23 no. 2 (2014), 483-512.

\end{thebibliography}
\end{document}